\documentclass[a4paper]{amsart}
  \usepackage{amssymb,exscale,bbm,mathrsfs}

 \newtheorem{theorem}{Theorem}[section]
 \newtheorem{corollary}[theorem]{Corollary}
 \newtheorem{lemma}[theorem]{Lemma}
 \newtheorem{proposition}[theorem]{Proposition}
 \theoremstyle{definition}
 \newtheorem{definition}[theorem]{Definition}
 \theoremstyle{remark}
 \newtheorem{remark}[theorem]{Remark}
 
 \newtheorem{examples}[theorem]{Examples}
\newtheorem*{bibnote}{Bibliographical notes}
 \numberwithin{equation}{section}

\begin{document}

  \renewcommand{\Re}{\text{Re}}
  \renewcommand{\Im}{\text{Im}}
  \newcommand {\sfrac}[2] { {{}^{#1}\!\!/\!{}_{#2}}} 
  \newcommand{\pihalbe}{\sfrac{\pi}2}
  \newcommand{\einhalb}{\sfrac{1}2}
  \newcommand{\Borel}{{\mathcal Bo}}
  \newcommand {\Poisson}{{\mathscr P}}
  \newcommand {\CC}{\mathbb C}
  \newcommand {\NN}{\mathbb N}
  \newcommand {\ZZ}{\mathbb Z}
  \newcommand {\RR}{\mathbb R}
  \newcommand {\cB}{\mathcal B}
  \newcommand {\cD}{\mathcal D}
  \newcommand {\cF}{\mathcal F}
  \newcommand {\cL}{\mathcal L}
  \newcommand {\eins} {\mathbbm 1}
  \newcommand {\al}{\alpha}
  \newcommand {\la}{\lambda}
  \newcommand {\eps}{\varepsilon}
  \newcommand {\si}{\sigma}
  \newcommand {\Si}{\Sigma}
  \newcommand {\Ga}{\Gamma}
  \newcommand {\ga}{\gamma}
  \newcommand {\om}{\omega}
  \newcommand {\Om}{\Omega}
  \newcommand {\Xm}{X_{-1}}
  \newcommand {\Sec}[1] {S({#1})}
  \newcommand {\range}{\textit{ ran }}
  \newcommand {\ran} {\textit{ ran }}
  \newcommand {\BOUNDED}{\mathcal B}
  \newcommand {\DOMAIN}{\mathcal D}
  \newcommand {\FOURIER}{\mathcal F}
  \newcommand {\LAPLACE}{\mathcal L}
\newcommand{\norm}[1] {\| #1 \|}
\newcommand{\lrnorm}[1]{\left\| #1 \right\|}
\newcommand{\bignorm}[1]{\bigl\| #1 \bigr\|}
\newcommand{\Bignorm}[1]{\Bigl\| #1 \Bigr\|}
\newcommand{\Biggnorm}[1]{\Biggl\| #1 \Biggr\|}
\newcommand{\biggnorm}[1]{\biggl\| #1 \biggr\|}
  \newcommand {\Bigidual}[3] {\Bigl\langle #1, #2 \Bigr\rangle_{#3}}
  \newcommand {\bigidual}[3] {\bigl\langle #1, #2 \bigr\rangle_{#3}}
  \newcommand {\idual}[3] {\langle #1, #2 \rangle_{#3}}
  \newcommand {\Bigdual}[2] {\Bigidual{#1}{#2}{}}
  \newcommand {\bigdual}[2] {\bigidual{#1}{#2}{}}
  \newcommand {\dual}[2] {\idual{#1}{#2}{} }
  \newcommand {\ws} {weak$^*$}
  \newcommand {\weg} {\backslash}
  \newcommand {\bbrace}[1] {[\![ #1 ]\!]} 
\def\eqnaintertext#1{\crcr\noalign{\vskip\belowdisplayskip
 \vbox{\normalbaselines\noindent#1}\vskip\abovedisplayskip}}
\newcommand {\suchthat}{:\;}
\newcommand {\emb} {\hookrightarrow}
\renewcommand{\labelenumi} {(\alph{enumi})}    
\renewcommand{\labelenumii}{(\roman{enumii})}
\renewcommand{\theenumi} {(\alph{enumi})}      
\renewcommand{\theenumii}{(\roman{enumii})}    
\allowdisplaybreaks
\title[On the Carleson measure criterion in linear systems theory]
 {On the Carleson measure criterion in linear systems theory}
\author[Bernhard H. Haak]{Bernhard H. Haak}
\address{
Institut de Math\'ematiques de Bordeaux\\Universit\'e Bordeaux 1\\351, cours de la Lib\'eration\\33405 Talence CEDEX\\FRANCE}
\email{bernhard.haak@math.u-bordeaux1.fr}
\date{\today}

\begin{abstract}
In Ho, Russell \cite{HoRussell}, and Weiss \cite{Weiss:Carleson}, a
Carleson measure criterion for admissibility of one-dimensional input
elements with respect to diagonal semigroups is given. We extend
their results from the Hilbert space situation
$X=\ell_2$ and $L^2$--admissibility to the more general situation of
$L^p$--admissibility on $\ell_q$--spaces.  In case of analytic
diagonal semigroups we present a new result that does not rely on
Laplace transform methods. A comparison of both criteria 
leads to result of $L^p$--admissibility for reciprocal systems 
in the sense of Curtain  \cite{Curtain:reciprocal}.
\end{abstract}
\maketitle

\section{Introduction}
Consider the infinite dimensional linear system described by the
differential equation
\begin{equation}
  \label{eq:control-system}\left\{
  \begin{array}{lcl}
    x'(t) + A x(t) &=& B u(t)\cr
    x(0) &=& x_0 \in X
  \end{array}\right.
\end{equation}
on a Banach space $X$. We assume that $-A$ generates a strongly continuous
semigroup $S(\cdot)$ on $X$. Then the spectrum of $A$ is necessarily
contained in some half plane $\CC_\beta$, $\beta\in\RR$. 
Choosing some element $\la$ of the resolvent set $\varrho(-A)$, we may
define the Banach space $X_1 = (\DOMAIN(A), \norm{(\la{+}A)\cdot})$, 
and the so-called extrapolation space $\Xm$ defined as the completion of $X$ 
with respect to a resolvent-norm $\norm{(\la{+}A)^{-1}\cdot}$ for 
some fixed $\la\in\varrho(-A)$  (see \cite{EngelNagel, Haase:Buch}
for more details on these spaces). Since $X$ is dense in $\Xm$, the 
semigroup $S$ extends in a natural way to $X_{-1}$; for the sake of
simplicity we denote this extension also $S$.

Let $U$ be another Banach space, and assume that $B\in
\BOUNDED(U,\Xm)$.
A solution of (\ref{eq:control-system}) is necessarily of the form
\[
x(t) = S(t) x_0 + \int_0^t S(t{-}s) B u(s)\,ds
\]
Notice that $x(t)$ is a well-defined element of $\Xm$ for $t\ge 0$ but 
that, in general, there is no reason why $x(t)$ should be an element of $X$. 

\begin{definition}
Let $p \in [1,\infty]$. A bounded operator $B\in \BOUNDED(U,\Xm)$ is called
{\em finite-time $L^p$--admissible} for $A$, if for
every $\tau>0$ there exists a constant $K>0$ such that for all $u \in
L^p([0,\tau], U)$
\begin{equation}\label{eq:finite-time-admiss}
\biggnorm{\int_0^t S(t-s) B u(s)\,ds}_X
\leq K
\norm{u}_{L^p(0,\tau; U)} \qquad t\in [0,\tau].
\end{equation}
The integral here is considered as a Bochner integral in $\Xm$ and we
suppose that the integral takes values in $X$.  The operator $B$ is
called {\em (infinite-time) $L^p$--admissible} for $A$, if
the constant $K$ in (\ref{eq:finite-time-admiss}) can be chosen
independently of $\tau>0$.
\end{definition}

In the case where the rank one
operator $B$, defined by $Bu(t) = b u(t)$ for $b\in \Xm$, is
(finite-time) $L^p$--admissible, we say that $b$ is an
$L^p$--admissible input {\em element}. 
When $p<\infty$, a density argument shows that $B$ is (infinite-time)
$L^p$--admissible if and only if the estimate
\[
\biggnorm{\int_0^\infty S(s) B u(s)\,ds}_X
\leq K
\norm{u}_{L^p(0,\infty; U)},
\]
holds for all $u$ in the Schwartz class. We remark that due to the
semigroup property, $L^p$--admissibility in finite and infinite time
coincide for uniformly exponentially stable semigroups. For the
special case $p=2$ there is a large literature on the notion of
admissibility. Among the early abstract formulations of admissibility
we refer e.g. to \cite{Salamon,Weiss:Admissibility-of-unbounded,%
Weiss:admiss-observation,Weiss:conjectures}.
A more recent publication \cite{JacobPartington:survey} gives a
detailed survey of the subject.

\medskip

In this article we restrict our focus to sequence spaces $X=\ell_q$.
Moreover, we assume $-A$ to be the generator of a bounded diagonal
semigroup $S(\cdot)$ where $\bigl( S(t)x\bigr)_n = \exp(-\la_n t)
x_n$, $n\in\NN$.  We consider the control system
\begin{equation}   \label{eq:control-system-eindim}
  x'(t) + A x(t) = b u(t), \qquad \text{for }t>0 \qquad \text{and} \quad x(0)=0
\end{equation}
with input element $b=(b_n) \in X_{-1} := \bigl\{ (\xi_n) \suchthat
\bigl(\frac{\xi_n}{1+\la_n}\bigr) \in X \bigr\}$.

Our aim is to give conditions for $L^p$-admissibility of $b$ for $A$
on $X=\ell_q$, i.e. conditions to guarantee the estimate

\[
   \biggl(\,\sum_{n=1}^\infty |b_n|^q \Bigl| \int_0^\infty e^{-\la_n s}
      u(s)\,ds\Bigr|^q\biggr)^\sfrac1q \le K \norm{u}_p
\]
or the respective estimate for convolutions derived from
(\ref{eq:finite-time-admiss}) with a uniform constant $K$ for all
times $t>0$. We present two criteria, both via a (generalised)
Carleson measure property of a discrete measure, associated with the
numbers $b_n$ and the eigenvalues $\la_n$ of $A$.

The first result is a generalisation of a result of {\sc Ho}, {\sc
  Russell} and {\sc Weiss}. It provides a sufficient condition for
$L^p$--admissibility in the case $p\le 2$ that characterises
$L^p$--admissibility if in addition $p' \le q$ where $p'$ is the
conjugate exponent defined by $\sfrac1p+\sfrac1{p'}=1$. The second
result is new and treats the case of analytic semigroups. It does not
rely on Laplace transform methods and therefore allows $p$ and $q$ to
be chosen freely in $(1,\infty)$.  It provides a sufficient condition
for all $p$ and $q$ that is necessary in the case $p\le q$. Notice
that, in contrast with the first result this allows us to characterise
$L^p$--admissibility on $X=\ell_p$ for all $p\in (1,\infty)$.

In the last section we discuss the two criteria. This leads naturally
to so-called \lq{}reciprocal systems\rq{} and the question of whether
admissibility of a system implies (or is implied by) admissibility of
the associated reciprocal system. We present two results, one for
general strongly continuous semigroups and the other for analytic
semigroups.

\section{Preliminaries on $\al$--Carleson measures}\label{sec:carleson}

Let $\RR^{d+1}_+ := \RR^d \times (0,\infty)$. Let 
$P_{t}(x) = c_d t (t^2+\norm{x}^2)^{-\frac{d+1}2}$
be the Poisson kernel on $\RR^d$ and let $(\Poisson f)(x,t)=\int
P_t(x-y)f(y)\,dy$ be the Poisson extension of $f$ to the half space
$\RR^{d+1}_+$. For $x \in \RR^d$, let $\Ga(x)$ denote the cone
$\{(y,t) \in \RR^{d{+}1}_+ \suchthat \norm{x-y} < t \}$. Given an open
set $O \subseteq \RR^d$ let $T(O)$ denote the ``tent''
\[
T(O) := \biggl( \bigcup_{x \not \in O} \Ga(x) \biggr)^\complement.
\]
If $d=1$ and $O$ is an open interval, this can be visualised as an
isosceles triangle with base $O$ in the half plane.  Let
$\Borel(\RR^{d+1}_+)$ be the set of all non-negative Borel-measures on
$\RR^{d+1}_+$.

\begin{definition}[embedding $\al$-Carleson measures]
Let $\al>0$ and $q \in [1,\infty)$ such that $\al q >1$. 
A non-negative measure $\mu \in \Borel(\RR^{d+1}_+)$ 
satisfying 
\begin{equation}   \label{eq:embedding-alpha-Carleson}
 \bignorm{f}_{L^q(\RR^{d+1}_+,\mu)} \le M_q \bignorm{f}_{H^{\al q}(\RR^{d+1}_+)} 
\end{equation}
for all $f \in H^{\al q}(\RR^{d+1}_+)$ is called an {\em embedding $\al$--Carleson  measure}.
\end{definition}

The above definition is independent of $q$; indeed, if 
$f \in L^{\al p}(\RR^d)$, then its Poisson extension $|P_t \ast
f|^{p/q}$ is subharmonic and therefore,  $|P_t \ast f|^{\sfrac{p}q}
\le P_t \ast |f|^{\sfrac{p}q}$.  This implies that
\[
    \bignorm{P_t \ast f}_{L^p(\mu)}^p 
\le \bignorm{ P_t \ast |f|^{\sfrac{p}q} }_{L^q(\mu)}^q
\le C \bignorm{ |f|^{\sfrac{p}q} }_{L^{\al q}(\RR^d)}^q 
 =  \bignorm{f}_{L^{\al p}}^p. 
\]

In case $\al=1$ we simply speak of Carleson measures. Recall that
$u \in H^q(\RR^{d+1}_+)$ if, and only if $u = \Poisson f$ for some 
$f\in L^q(\RR^d)$ (see, e.g. \cite[III.4.2]{Stein:singular-int})
and that $\norm{\Poisson f}_{H^q} = \norm{f}_{L^q}$.

\begin{definition}[geometric $\al$-Carleson measures]
A non-negative Borel measure $\mu \in \Borel(\RR^{d+1}_+)$ is  
called {\em geometric $\al$--Carleson} if
\begin{equation}   \label{eq:geometric-al-Carleson}
 \mu( T(Q) )^\al \leq c \,|Q|
\end{equation}
for all cubes $Q \subseteq \RR^d$.
\end{definition}

In the literature the choice of the exponent is not consistent
(sometimes $\al$ is replaced by $\sfrac1\al$).  For all $\al>0$,
embedding $\al$--Carleson measures are also geometric $\al$--Carleson;
this can be seen by applying the embedding estimate to the function
$P_t \ast f$ where $f=\eins_Q$ is the characteristic function of $Q$.
A celebrated theorem of {\sc Carleson} \cite{Carleson:interpolation,
  Carleson:corona} states that both notions coincide for $\al=1$, {\sc
  Duren} \cite{Duren:Carleson} extended this result to $\al\in (0,1]$.
In case $\al>1$, the embedding $\al$--Carleson property is strictly
stronger than geometric one as {\sc Taylor} and {\sc Williams}
\cite{TaylorWilliams:interpolation} show with a counterexample on the
complex unit disk. In Example~\ref{ex:examples}~\ref{item:ex_b} below
we give a similar example on the half space $\RR^2_+$.  Both notions
were subsequently treated in {\sc Amar} and {\sc Bonami}
\cite{AmarBonami}. Variants of these results with weighted Hardy
spaces have been obtained by {\sc McPhail}
\cite{McPhail:interpolation} and {\sc Nakazi}
\cite{Nakazi:interpolation}. Applications of Carleson measures in the
context of systems theory can be found e.g. in recent publications in
\cite{JacobPartingtonPott:tangential} and      
\cite{Wynn:discrete-continuous-time}.

\smallskip

Let $L^{\al,1}(\RR^{d+1}_+)$ denote the smallest Lorentz space of $\al$-integrable 
functions. Then $L^{\al,1}(\RR^{d+1}_+) \subsetneq L^{\al,\al}(\RR^{d+1}_+) = L^\al(\RR^{d+1}_+)$ for $\al>1$. 
If, in place of (\ref{eq:embedding-alpha-Carleson}), one merely knows that
$\norm{ \Poisson f }_{L^1(\RR^{d+1}_+, \mu)} \le C \norm{f}_{ L^{\al,1}(\RR^{d+1}_+)}$,
$\mu$ is called {\em weakly embedding $\al$--Carleson}. We shall not go into details 
about this notion, but just mention that it is linked to the distinction 
between requiring the estimate (\ref{eq:geometric-al-Carleson}) on 
all cubes or on all open sets (see \cite[Theorem 1]{AmarBonami}), and
consequently, the weak embedding notion implies the geometric one. 
In case $\al \le 1$ all notions of $(\al)$--Carleson coincide by the
Carleson-Duren result; all three are strictly different in case $\al>1$. 

\medskip

In the case $d=1$,{\sc  Vidensi\u{i}} \cite{Videnskii} gives several
equivalent conditions for $\mu$ to be embedding $\al$--Carleson for
$\al>1$, one of which is the Fefferman-Stein maximal function
description. It follows from a  maximal inequality due to 
{\sc  Fefferman} and {\sc Stein} 
\cite{FeffermanStein:maximal-inequalities}. Indeed, let 
\begin{equation}\label{eq:fefferman-stein-function}
 \psi_\mu(x) := \sup\Bigl\{ \frac{ \mu( T(Q) ) }{|Q|} \suchthat x\in Q \Bigr\}
\end{equation}
where the sets $Q$ in the above supremum are cubes in $\RR^d$
containing the element $x$. We call $\psi_\mu$ the Fefferman-Stein maximal function.  
\begin{theorem}\label{thm:fefferman-stein}
Let $\al>1$ and let $\beta$ be its conjugate exponent satisfying
$\sfrac1\al+\sfrac1\beta =1$. Then $\mu \in \Borel(\RR^{d+1}_+)$ is embedding
$\al$--Carleson if and only if $\psi_\mu \in L^\beta(\RR^d)$.
\end{theorem}
\begin{proof}
Suppose that $\psi_\mu \in L^\beta(\RR^d)$. 
The Fefferman and Stein inequality \cite[Theorem 2]{FeffermanStein:maximal-inequalities}
reads: $\norm{ \Poisson f }_{L^p(\mu)} \le \norm{f}_{L^p(\psi_\mu(x)\,dx)}$.
Since $\psi_\mu \in L^\beta$, H\"older's inequality now implies that $\mu$
is embedding $\al$-Carleson.\\
Conversely, if $\mu$ is embedding $\al$-Carleson, $P_t\ast:  L^\al \to
L^1(\mu)$  is bounded, which is equivalent to $F \in L^\beta(\RR^d)$
where $F(y) := \int P_{t}(x-y) \,d\mu(x,t)$ is the \lq{}balay\'ee\rq{}
of $\mu$. Notice that $(x,t) \in T(Q)$ implies that $B(x, t) \subseteq Q$ and so
\begin{equation}  \label{eq:integrierter-Poisson-kern}
   \int_{Q}  P_{t} (x-y) \,dy \ge 
   \int_{B(x,t)} P_t(x-y)\,dy = 
   \int_{B(0,1)} P_1(y)\,dy =: \eps_d > 0
\end{equation}
for $(x,t) \in T(Q)$. Consequently, 
\[
    \frac1{|Q|} \int_Q F(y)\,dy 
 =  \frac1{|Q|} \int_{T(Q)} \int_Q P_t(x-y)\,dy \,d\mu(x,t) 
\ge \eps_d \frac{\mu(T(Q)}{|Q|}.
\]
Denoting by $M F$ the (uncentred) maximal function of $F$ the above
estimates yields $(M F)(x) \ge \eps_d \,\psi_\mu(x)$ and therefore
$\psi_\mu \in L^\beta(\RR^d)$.
\end{proof}

Another way of verifying that  $\mu$ is embedding $\al$--Carleson is to 
show that (\ref{eq:embedding-alpha-Carleson}) holds for linear combinations of 
reproducing kernels by a density argument. Indeed, a continuous linear
functional $l\in \bigl(H^p(\RR^d)\bigr)' = H^{p'}(\RR^d)$ 
that vanishes on all reproducing kernels $k_x$ satisfies
$0=l(k_x)=l(x)$ and therefore $l=0$. It is natural to ask whether 
the embedding $\al$--Carleson property can be tested on 
the reproducing kernels {\em  without} taking linear combinations. 
This is sometimes referred to as the \lq{}reproducing kernel
thesis\lq{}. 

\begin{lemma}\label{lem:geometric-RKT}
Assume that for $p,q \in (1,\infty)$ there exists a constant $M>0$
such that
\begin{equation}   \label{eq:kern-abschaetzung}
    \norm{ k_z }_{L^q(\RR^2_+),\mu)} \le M \norm{ k_z }_{H^p(\RR^2_+)}
\end{equation}
for all $z\in \CC_+$. Then $\mu$ is geometric $\al$--Carleson where
$\al = \sfrac{p}q$.
\end{lemma}
\begin{proof}
Since $k_\la(z) = 1/(z+\bar \la)$,
letting $\Re(\la) = \xi>0$ the substitution $y=\xi t$ shows
\[    \norm{ k_\la }^{p}_{H^{p}} 
= \int_{-\infty}^\infty \frac{dy}{(y^2+\xi^2)^{\sfrac{p}2}} \\
= \int_{-\infty}^\infty \frac{\xi \, dt}{\xi^{p}(t^2+1)^{\sfrac{p}2}}
 =   C_{p}^{p} \; \xi^{1-p},
\]
and therefore, $\norm{ k_\la }_{H^p} = C_p
(\Re(\la))^{-\sfrac1{p'}}$. For $\om \in \RR$ and $r>0$, let 
$\la = i\om + r$ and use the shorter notation $T_{\om, r}$ for
the tent $T( (\om{-}r, \om{+}r))$.
Then the triangle inequality yields $|k_\la(z)| \ge 1/(2r)$ for all $z
\in T(\om{-}r, \om{+}r)$ and therefore, 
\begin{align*}
  \mu( T_{\om, r} ) 
&=   \int_{ T_{\om, r} } \,d\mu 
\le (2r)^q \int_{ T_{\om, r} } \bigl| k_\la(z)\bigr|^q \,d\mu \\
&\le M (2r)^q \bignorm{ k_\la(z) }_{H^p}^q 
 =  M C_p 2^q  r^{q-\sfrac{q}{p'}} = C \; r^{\sfrac1\al}
\end{align*}
This shows that $\mu$ is geometric $\al$--Carleson as claimed. 
\end{proof}

Using the Lemma and the results from Duren-Carleson and the
aforementioned counterexample from Taylor and Williams, it becomes
clear that a $H^p$-$L^q(\mu)$ version of the reproducing kernel thesis
holds if, and only if  $p \le q$.

\begin{examples} \label{ex:examples} 
Let $d=1$ and identify $\RR^2_+$ with the right half plane $\CC_+$.\hfill
\begin{enumerate}
\item \label{item:ex_a}
     The Lebesgue measure in $\RR^2_+$ is clearly geometric
     $\einhalb$--Carleson and therefore also embedding $\einhalb$--Carleson.
\item \label{item:ex_d}
     For $\om \in(0,\pi)$, let  $S(\om) := \{z\in\CC
     \setminus\{0\}:|\arg z| < \om\}$ and let $S(0) := (0,\infty)$.
      We say that $\mu$ is sectorial, if the support of $\mu$ is contained 
      in a finite union of sectors $i \om_k+S(\theta_k)$, $k=1\ldots,n$ with 
      $\om_k\in \RR$ and $\theta\in (0,\pihalbe)$. \\
      If $\mu$ is sectorial and geometric $\al$--Carleson, intervals
      $I$ containing a point $x$ with sufficiently large absolute value need
      a length of at least $\eps |x|$ for $T(I)$ to intersect with the
      support of $\mu$. Therefore, geometric $\al$--Carleson implies
      that the Fefferman-Stein maximal function is weak--$L^\beta$.
      Consequently, sectorial measures $\mu$ that are geometric $\al$-- 
      and geometric $\widetilde \al$--Carleson are embedding $\ga$--Carleson for
      all $\ga \in (\al, \widetilde \al)$ by Marcinkiewicz'
      interpolation theorem.  
\item \label{item:ex_e}
      Sectorial measures $\mu$  are embedding
      $\al$--Carleson if  and only if the function $g$ defined by 
       $g(r) := \frac1r \,\mu(\{ \Re(z) < r\})$ satisfies $g \in L^\beta(\RR_+)$
      (see \cite[Remark, p. 188]{Videnskii}). 
\item \label{item:ex_b}
     For  $\al>1$, measures that are geometric- but not embedding
     $\al$--Carleson  are given in \cite{TaylorWilliams:interpolation}
     for the unit disc.  A simplified half-space version reads as
     follows: let 
     $\al\in (1,2)$ and let $\mu_\al$ be the measure
      $\mu_\al = \bigl(\frac{x}{1+x}\bigr)^{\sfrac1\al-1} \,dx$,
      supported on  the real line. Then, for every interval $I$, 
      one has $\mu_\al(T(I)) \le C\, |I|^\sfrac1\al$, 
      whence $\mu_\al$ is geometric $\al$--Carleson. 
      For a real parameter $r>0$ let $F(z) = (r+z)^{-2/p}$. 
      For $p,q \in (1,\infty)$   with $\al=\sfrac{p}q$ one has 
     \[
       \norm{F}_{H^p} 
      = \biggl( \int_{-\infty}^\infty \frac1{r^2+t^2} \,dt \biggr)^\sfrac1p
      = c_p \,r^{-\sfrac1p},
     \]
      whereas  a  substitution yields
     \[
       \norm{F}_{L^q(\mu_\al)} 
      = r^{-\sfrac1q} \biggl( \int_{0}^\infty \bigl(\frac{t}{1+rt}\bigr)^{\sfrac{q}p-1}
                     \frac1{(1+t)^{2q/p}} \,dt \biggr)^\sfrac1q .
     \]
     By  Lebesgue's dominated convergence theorem, one has 
     $\norm{F}_{L^q(\mu_\al)}  \sim r^{-\sfrac1q}$ for $r\to
     0+$. Therefore, the embedding $\al$--Carleson property fails.
\item \label{item:ex_c}
     Finally we give an example of a {\em discrete} measure that is
     geometric $\al$--Carleson 
     without being embedding $\al$--Carleson: let $\eps \in [0, 1)$ and
     $\ga\ge 1$ such that $\al=\tfrac\ga{1-\eps} > 1$. Let $\mu =
     \sum_{n  \in\ZZ^*} |n|^{-\eps} \delta_{\la_n}$ where $\la_n = 1+i n^\gamma$. 
     Since $\mu(T(B(x, r)))=0$ for $r<1$ we may suppose $r\ge 1$.
     Since $\ga\ge 1$ the support of $\mu$ is either equidistant or
     thins out at infinity. Moreover weights decrease as $|n|$ increases. Therefore
\begin{align*}
     \mu\bigl(T( B(x, r))\bigr) 
&\le \mu\bigl(T(B(0, 2r))\bigr) 
 \le  2\sum_{j=0}^{\lfloor (2r{-}1)^{\sfrac1\ga} \rfloor}   j^{-\eps} \\
&\le 2 \int_1^{(2r)^{\sfrac1\ga}} (t-1)^{-\eps}  \,dt
 \le \frac2{1-\eps} (2r)^{\frac{1-\eps}\ga},
\end{align*}    
     whence $\mu$ is geometric $\al$--Carleson.
     Considering the tents
     $T((-x{-}1,x{+}1))$ for $x>1,$ however, yields the estimate 
\[
    \psi_\mu(x) \ge c \, \frac{(\lfloor  (x)^\frac1\ga
      \rfloor^{1-\eps}-1)}{x}. 
\]
Let $\beta$ denote the  dual exponent of $\al$. Then, for $n\ge 2$, it
readily follows that
\[
    \int_{n^\ga}^{(n{+}1)^\ga} \bigl|\psi_\mu(s)\bigr|^\beta \,ds \ge
    \frac{ c_{\eps,\beta,\ga} }{n+1}
\]
where $c_{\eps,\beta,\ga}>0$.
Now Theorem~\ref{thm:fefferman-stein} allows us to conclude that 
$\mu$ cannot be embedding $\al$--Carleson since the above estimate
shows $\psi_\mu \not\in L^\beta(\RR)$.
\end{enumerate}
\end{examples}

In order to apply Carleson-type criteria for admissibility,
we wish to replace the Poisson kernel by more general 
convolution kernels. In case $\al \le 1$ the
following result is due to {\sc Stein} \cite[Theorem
II.5.9]{Stein:harmonic-analysis}.

\begin{theorem}\label{thm:Lq-Lp-Carlseon}
Let $q \in (1, \infty)$ and $\al > \sfrac1q$.
Let $\Phi$ be a function on $\RR^d$ that admits a radial, non-increasing,
bounded and integrable dominating function $\varphi$. For $t>0$ put
$\Phi_t(x) := t^{-d} \Phi(x/t)$. Then for all embedding $\al$--Carleson measures
$\mu$ on $\RR^{d+1}_+$ the estimate
\[
\biggl(\int_{\RR^{d+1}_+} \bigl|(\Phi_t \ast f)(x)\bigr|^q
    \,d\mu(t,x)\biggr)^{\sfrac{1}q}
\leq c \norm{f}_{L^{\al q}(\RR^d)}
\]
holds for all $f\in L^{\al q}(\RR^d)$.
\end{theorem}
\begin{proof}
Notice that if $F(x, t) := (\Phi_t \ast f)(x)$ and  
$F_\Phi^*(x) := \sup_{|x-y|<t} |F(y, t)|$, then $s>0$
\[
  \{ (y, t)\in \RR^{d+1}_+ \suchthat |F(y,t)| > s \} \subseteq T(O_s)
\]
where $O_s :=\{ x \in \RR^d \suchthat F_\Phi^*(x) > s \}$, see, e.g., 
\cite[Section II.2.3]{Stein:harmonic-analysis}. Now
inequality~(\ref{eq:integrierter-Poisson-kern}) implies
$(\Poisson F_\Phi^*)(x,t) \ge s \, \eps_d$ for $(x,t) \in T(O_s)$. 
The problem now boils down to the Poisson kernel estimate as follows:
\begin{eqnarray*}
      \norm{ F }_{L^q(\mu)}^q 
&=&   q \int_0^\infty s^{q-1} \mu\left( |F|>s \right)\,ds 
\le   q \int_0^\infty s^{q-1} \mu\left( T(O_s) \right)\,ds \\
&\le& q \int_0^\infty s^{q-1} \mu\left( (\Poisson F_\Phi^*) > s \eps_d \right)\,ds 
 = \eps_d^{-q} \norm{ \Poisson F_\Phi^* }_{L^q(\mu)}^q.
\end{eqnarray*} 
The desired estimate 
\[
\norm{ F }_{L^q(\mu)} \le \tfrac{1}{\eps_d} \norm{\Poisson F_\Phi^*}_{L^q(\mu)} 
\le \tfrac{C}{\eps_d} \norm{F_\Phi^*}_{L^{\al q}(\RR^d)} \le \tfrac{C M}{\eps_d} \norm{f}_{L^{\al q}(\RR^d)}
\]
follows from Theorem~\ref{thm:fefferman-stein} and 
\cite[Proposition II.2.1]{Stein:harmonic-analysis}.
\end{proof}

\section{Criteria for admissibility of diagonal systems}\label{sec:criteria}
As a first result we present a Carleson measure criterion for
admissibility of diagonal systems. 

\begin{theorem}\label{thm:alpha-Carleson-Poisson}
Let $q\in (1, \infty)$, $p \in (1, 2]$ and $\al q = p'$ where 
$p'$ is the dual exponent of $p$. On $X=\ell_q$ let $A$ be a
diagonal operator with eigenvalues $\la_n \in \CC_+$ and let $b=(b_n)$ a
sequence of complex numbers.
Consider the discrete measure  $\mu = \sum_n |b_n|^q
\delta_{\la_n}$.
\begin{enumerate}
\item\label{item:alpha-Carleson-Poisson-a} 
       If $\mu$ is an embedding $\al$--Carleson measure, then 
       $b \in X_{-\theta}$ for all $\theta > \sfrac1{p'}$ and $b$ is
       an infinite-time $L^p$--admissible input element for $A$.
\item \label{item:alpha-Carleson-Poisson-b} 
       If $b$ in an infinite-time $L^p$--admissible input element for
       $A$, then $b \in X_{-\theta}$ for all $\theta>\sfrac1{p'}$ and 
       $\mu$ is geometric $\al$--Carleson.
\end{enumerate}
In particular, the $\al$--Carleson property of $\mu$ characterises the
$L^p$--admissibility of $b$ in the case $p'\le q$. Moreover, the
result in \ref{item:alpha-Carleson-Poisson-b} is optimal in the sense
that for $p'>q$, the measure $\mu$ is not embedding $\al$--Carleson in
general.
\end{theorem}
\begin{proof}
\ref{item:alpha-Carleson-Poisson-a}
  First assume that $\mu$ is an embedding $\al$--Carleson measure.
  Then $\mu$ is also geometric $\al$-Carleson. To show $b \in
  X_{-\theta}$ we employ essentially the argument in \cite[Proposition
  3.1]{Weiss:Carleson}: let $\Delta_n = T\bigl( [{-}2^n, 2^n] \bigr)
  \setminus T\bigl( [{-}2^{n-1}, 2^{n-1}] \bigr)$. Then
\begin{align*}
\sum_k \frac{ |b_k|^q }{|1+\la_k|^{\theta q} }
& \le \sum_{n\in \ZZ} \frac1{(1+2^{n-1})^{\theta q}} \sum_{\la_k \in
  \Delta_n} |b_k|^q \\
& \le \sum_{n\in \ZZ} \frac1{(1+2^{n-1})^{\theta q}} \;\mu\bigl( T([{-}2^{n},
2^n]) \bigr)
\le C \sum_{n\in \ZZ} \frac{2^{\frac{nq}{p'}}}{1+2^{n \theta q}}
\end{align*}
which is finite for all $\theta>\sfrac1{p'}$. 

To show that $b$ is infinite-time $L^p$--admissible, let $\tau>0$ and
let $u$ be a test function with support in $[0,\tau]$. Then its
Laplace transform $\LAPLACE u$ is bounded and analytic on the right half
plane and therefore can be reproduced by its (non-tangential) boundary
values on $i \RR$. Notice that $(\LAPLACE u)(is) = \FOURIER u(s)$ where $\FOURIER u$
denotes the Fourier transform of $u$.
\begin{eqnarray*}
     \biggnorm{\int_0^\infty S(s) b u(s) \,ds}_{\ell_q}
&=&  \biggl(\sum_{n=0}^\infty |b_n|^q \biggl| \int_0^\infty
       e^{-\la_n s} u(s) \, ds\biggr|^q\biggr)^{\sfrac{1}q}\\
(\text{definition of }\mu)
&=&  \biggl( \int_{\RR^2_+} \bigl| (\LAPLACE u)(t+ix) \bigr|^q
      \,d\mu(x,t)\biggr)^{\sfrac{1}q} \\
&=&  \bignorm{ P_t\ast (\FOURIER u) }_{L^q(\mu)} \leq  C \, \bignorm{ \FOURIER u }_{L^{p'}(\RR)} \\
&\leq&  C \, \norm{ u }_{L^p(\RR_+)}.
\end{eqnarray*}
In the last estimate we make use of $p\le 2$ and the boundedness of
the Fourier transform from $L^p$ to $L^{p'}$. Notice that the obtained
constant is independent of $\tau>0$.

\medskip

\ref{item:alpha-Carleson-Poisson-b} Conversely assume $b$ to be $L^p$--admissible for some
$p\in(1,\infty)$ and consider the reproducing  kernel function
$k_z(\la) = \frac1{\la+\overline z}$ for $\Re(z)>0 $. Then
\begin{equation}\label{eq:embedding}
  \begin{split}
      \bignorm{ k_z }_{L^q(\mu)} 
& =  \biggl( \sum_{n=0}^\infty \left|\frac{b_n}{\la_n+{\overline{z}}} \right|^q \biggr)^\sfrac1q
  =  \biggnorm{\int_0^\infty T(t) b e^{-t \,\overline{z}} \,dt}_{\ell_q} \\
&\le C \bignorm{ e^{-\cdot \,\overline{z}} }_{L^p(\RR_+)}
  =  C \biggl( \Re(z) \biggr)^{-\sfrac1p}
  =  C'\,  \norm{k_z}_{H^{p'}(\CC_+)}.
  \end{split}
\end{equation}
So, by Lemma~\ref{lem:geometric-RKT}, $\mu$ is geometric
$\al$--Carleson. Consequently, $b \in X_{-\theta}$ for all
 $\theta>\sfrac1{p'}$.

\smallskip

It remains to show that $\mu$ is not embedding $\al$--Carleson in
general. We use Example~\ref{ex:examples}~\ref{item:ex_c}. 
For a given $q \in(1,\infty)$ we choose $\eps\in (0,1)$ such that
$q\,\eps > 1$. Then $b=(b_n)$ with $b_n = n^{-\eps}$ will satisfy $b\in
\ell_q$. For a  given $\al>1$, let $\ga$ and $(\la_n)$
be chosen as in the example. Let $p' = \al q$. Using 
H\"older's inequality and the fact that $\Re(\la_n)=1$,
\begin{eqnarray*}
  \biggnorm{ \int_0^\infty S(t)b u(t) \,dt}_{\ell_q}
&=& \biggl( \sum_{n=1}^\infty \bigl|b_n\bigr|^q \biggl| 
    \int_0^\infty e^{-\la_n t}   u(t)\,dt \biggr|^q\biggr)^\sfrac1q \\
&\le& \norm{u}_{L^p(\RR_+)} \; \biggl(\sum_{n=1}^\infty \bigl|b_n\bigr|^q
      \norm{e^{-t}}_{L^{p'}(\RR_+)}^q \biggr)^\sfrac1q,
\end{eqnarray*}
and so $b$ is $L^p$ admissible. 
However, as shown in
the Example~\ref{ex:examples}~\ref{item:ex_c}, the embedding
$\al$--Carleson property fails. 
\end{proof}

\medskip

\begin{bibnote}
Part \ref{item:alpha-Carleson-Poisson-a} of the theorem in
the case $p=q=2$ is due to {\sc  Ho}, {\sc Russell} 
\cite{HoRussell} and the case $\al\le1$, in which the geometric and
the embedding $\al$--Carleson properties coincide, was found by 
{\sc  Unteregge} \cite{Unteregge} independently.
Part \ref{item:alpha-Carleson-Poisson-b} of the theorem is due to {\sc Weiss}
\cite{Weiss:Carleson}, but we provide a different proof. 
\end{bibnote}

\smallskip

In the proof of Theorem~~\ref{thm:alpha-Carleson-Poisson}, the
embedding $\al$--Carleson property  is 
applied somehow indirectly: first the desired estimate is transformed into 
an interpolation problem for the Laplace transform of $u$, which makes
it necessary to use boundedness of the (inverse) Fourier
transform; in the end the unaesthetic restriction $p\le2$ seems inevitable.
This is particularly embarrassing since $p=q$ is a very natural case
to consider.

\smallskip

The next result overcomes this obstacle; it makes no use of the
Laplace transform but  instead directly uses the convolution estimate from 
Theorem~\ref{thm:Lq-Lp-Carlseon}. However, it requires analyticity of the
semigroup $S(\cdot)$ as an additional assumption. 
We recall some notation: we recall that a  
\emph{sectorial operator $A$ of type $\om\in[0,\pi)$} in a Banach
space $X$ is a closed linear operator $A$ satisfying 
$\sigma(A)\subseteq \overline{S(\om)}$ and, for any $\nu\in(\om,\pi)$,
\[
   \sup \bigl\{ \,\|\la R(\la,A)\|: |\arg\la|\ge\nu \bigr\} < \infty.
\]
An operator $-A$ generates a bounded analytic semigroup in
$X$ if and only if  $A$ is a densely defined sectorial operator in $X$
of type $<\pihalbe$. 
\begin{theorem}\label{thm:alpha-Carleson-direkt}
Let $p,q \in (1, \infty)$ and $\al q = p$. Let $\theta\in
(0, \pihalbe)$ and let $A$ be an injective diagonal operator on
$X=\ell_q$ with eigenvalues $\la_n \in S(\theta)$. Let
$b = (b_n)$ be a sequence of complex numbers and consider the discrete
measure  $\mu  = \sum_n |\tfrac{b_n}{\la_n}|^q \delta_{\lambda_n^{-1}}$.
\begin{enumerate}
\item \label{item:alpha-Carleson-direkt-a} 
  If $\mu$ is embedding $\al$--Carleson, then $b \in X_{-\theta}$ for
  all $\theta>\sfrac1{p'}$ and $b$ is an (infinite-time)
  $L^p$--admissible input element for $A$.
\item \label{item:alpha-Carleson-direkt-b}
  If $b$ is infinite-time $L^p$--admissible for $A$, then $b \in
  X_{-\theta}$ for all $\theta>\sfrac1{p'}$ and $\mu$ is geometric
  $\al$--Carleson.
\end{enumerate}
In particular, the $\al$--Carleson property of $\mu$ characterises
$L^p$--admissibility of $b$ in the case $p \le q$.
\end{theorem}
\noindent Notice that if $0\in\varrho(A)$, then $\bigl(\tfrac{b_n}{\la_n}\bigr)_n \in \ell_q$
and consequently  $\mu$ is a finite measure.
\begin{proof}
  \ref{item:alpha-Carleson-direkt-a} Let $\mu$ be an embedding
  $\al$--Carleson measure.  To see that $b \in X_{-\theta}$ for all
  $\theta>\sfrac1{p'}$ is essentially the same as in the proof of
  Theorem~\ref{thm:alpha-Carleson-Poisson} and so we omit it. 
  To show infinite-time admissibility, choose an interval $[0,\tau]$
  with $\tau>0$ and let $u\in L^p(0,\tau)$. We shall prove
  (\ref{eq:finite-time-admiss}) with a constant $K$ that is
  independent of $\tau>0$.

  Let $r_n = \Re(\la_n)$ and let $\widetilde \mu = \sum_n
  \bigl|\tfrac{b_n}{\la_n}\bigr|^q \delta_{(r_n^{-1}+i\tau)}$. We
  claim that $\widetilde \mu$ is an embedding $\al$--Carleson-measure
  on $\RR^2_+$ as well. This has nothing to do with $\tau$ and so it
  is sufficient to show that $\mu$ is an embedding $\al$--Carleson
  measure on $\RR^2_+$ if and only if $\nu = \sum_n
  \bigl|\tfrac{b_n}{r_n}\bigr|^q \delta_{r_n^{-1}}$ is an embedding
  $\al$--Carleson measure. To prove this claim we use the shorter
  notation $T_{\om,r}$ for the tent $T\bigl( (\om{-}r,\om{+}r)\bigr)$.
  A simple geometric consideration, along with the fact that on the
  sector
  $S(\theta)$, %
  real parts and absolute values are equivalent up to a constant of
  $\frac1{\cos(\theta)}$, shows that
  \[
  \mu\bigl( T_{w, r} \bigr) \le \nu\bigl( T_{0, r} \bigr) \le
  \tfrac1{\cos(\theta)^q} \; \mu\bigl( T_{0, |w|+r} \bigr).
  \]
  If $\al \le 1$ and if $\nu$ is $\al$--Carleson, the above inequality
  immediately yields that $\mu$ is $\al$--Carleson. Next, suppose that
  $\al\le1$ and that $\mu$ is $\al$--Carleson.  Observe, that we may
  suppose for any tent $T_{\om,r}$ that $|\om| \le \tan(\theta)\,r$
  since otherwise there is no intersection of the tent with the sector
  $S(\theta)$. This implies the estimate
\[
    \nu\bigl( T_{\om, r} \bigr) 
\le \nu\bigl( T_{0, r}\bigr) 
\le \tfrac1{\cos(\theta)^q} \; \mu\bigl( T_{0, |w|+r} \bigr) 
\le \tfrac{C}{\cos(\theta)^q} \; (|\om|+r)^\sfrac1\al 
\le  \tfrac{C(1{+}\tan(\theta))}{\cos(\theta)^q} \; r^\sfrac1\al
\]
  and so $\nu$ is $\al$--Carleson.  Finally, if $\al>1$, the measures
  $\mu$ and $\nu$ are simultaneously embedding $\al$--Carleson by
  Example~\ref{ex:examples}~\ref{item:ex_e}. This proves the claim.

\smallskip

We now know that $\widetilde \mu$ is an embedding $\al$--Carleson
measure.  Consider the kernel function
\[
  \Phi(x) := \exp(-x) \eins_{[0,\infty)}(x)
\]
which satisfies the assumptions of Theorem~\ref{thm:Lq-Lp-Carlseon} (a
radially decreasing integrable majorant is $\exp(-|x|)$). Notice that
for $s>0$, $\Phi_s(x) = s^{-1} \exp(-x/s) \eins_{[0,\infty)}(x)$. Then
\begin{eqnarray*}
     \biggnorm{\int_0^\tau S(s) b u(\tau{-}s) \,ds}_{\ell_q}
& = &\biggl(\sum_{n=0}^\infty \bigl|\tfrac{b_n}{\la_n}\bigr|^q
     \biggl| \int_0^\tau \la_n e^{-\la_n s} u(\tau{-}s) \, ds\biggr|^q\biggr)^{\sfrac{1}q}\\
&\leq&\biggl(\sum_{n=0}^\infty \bigl|\tfrac{b_n}{\la_n}\bigr|^q
     \biggl( \int_0^\tau \bigl|\tfrac{\la_n}{r_n}\bigr| r_n e^{-r_n s}
        |u(\tau{-}s)| \, ds\biggr)^q\biggr)^{\sfrac{1}q}\\
&\leq& \tfrac1{\cos(\theta)}
     \biggl(\sum_{n=0}^\infty \bigl|\tfrac{b_n}{\la_n}\bigr|^q
       \bigl((\Phi_{r_n^{-1}} \ast |u|)(\tau)\bigr)^q\biggr)^{\sfrac{1}q}\\
(\text{definition of }\widetilde \mu)
& = &  \tfrac1{\cos(\theta)}
      \biggl(\int_{\RR^2_+} \bigl( (\Phi_s \ast |u|)(t) \bigr)^q
      \,d\widetilde \mu(s,t)\biggr)^{\sfrac{1}q}\\
& \leq &\tfrac{c}{\cos(\theta)} \norm{u}_{L^p}
\end{eqnarray*}
where we used Theorem~\ref{thm:Lq-Lp-Carlseon} in the last estimate.
Since the estimate is independent of the choice of $\tau>0$, $b$ is
infinite-time admissible for $A$.

\smallskip

\ref{item:alpha-Carleson-direkt-b} Now assume that $b$ is
$L^p$--admissible. Then, for $\Re(z)>0$, the reproducing kernel
functions $k_z(\la) = \frac1{\la+\overline z}$ satisfy
\begin{eqnarray*}
     \bignorm{ k_z }_{L^q(\mu)} 
& = &\biggl( \sum_{n=0}^\infty \left|\frac{b_n}{\la_n} 
        \frac1{\la_n^{-1}+{\overline{z}}} \right|^q \biggr)^\sfrac1q
  =  \frac1{|z|} \biggl( \sum_{n=0}^\infty \left| \frac{b_n}{\la_n+\overline{z}^{-1}} 
         \right|^q \biggr)^\sfrac1q \\ 
& = & \frac1{|z|} \biggnorm{\int_0^\infty T(t) b  e^{-t \overline{z}^{-1} } \,dt}_{\ell_q} 
 \le  \frac{C}{|z|} \bignorm{ \exp(-t \overline{z}^{-1} )}_{L^p(\RR_+)}\\
& = & \frac{C}{|z|} \biggl( \frac{|z|^2}{\Re(z)} \biggr)^\sfrac1p 
 \le  C \, \cos(\theta)^{-\sfrac2p} \, |\Re(z)|^{-\sfrac1{p'}}\\
& = & C \, \cos(\theta)^{-\sfrac2p} \; \norm{k_z}_{H^p(\RR^2_+)},
\end{eqnarray*}
and so Lemma~\ref{lem:geometric-RKT} yields that $\mu$ is geometric $\al$--Carleson.
\end{proof}

\section{$L^p$--admissible control operators}

We now consider the case that $U=X=\ell_q$. Let $B: U \to X_{-1}$ be
linear and bounded. 
Then there are functionals $\varphi_n \in (\ell_q)^\ast$
such that $(B u)_n = \dual{\varphi_n}{u}$. Indeed, $B$
is determined by its values on any basis $(e_j)$ of $U$. We choose
$(e_n)$ to be the standard basis in $\ell_q$.  Let $\varphi_n$ be the
sequence $\bigl(\dual{ e_n} {B e_j}\bigr)_{j=1}^\infty$,
i.e. $\varphi_n$ is the scalar sequence of $n^{\rm th}$ coordinates of the
vector-valued sequence $(B e_j)$.  Then, for  $u = (u_j) \in U$, and $n\in \NN$,
\[
\dual{ e_n}{ B u } = \Bigdual{ e_n}{ \sum u_j B e_j } = 
\sum u_j \bigdual{ e_n}{  B e_j }
\]
is a finite number whence $\varphi_n \in (\ell_q)^\ast$ for all $n$.
So, $\dual{ e_n} { B u } = \dual{\varphi_n}{u}$, which proves the
claim. This means that any bounded linear operator $B: U \to X_{-1}$ may
be interpreted as a certain sequence of functionals. The following
proposition is a direct generalisation of \cite[Proposition
4.8.6]{TucsnakWeiss}.

\begin{proposition}\label{prop:extension-to-infinite-dim-U}
  Let $X = U = \ell_q$ and $1 < p \le q <\infty$. Let $(\varphi_n)$ be a
  sequence of elements in $U^\ast = \ell_{q'}$ and consider the
  scalar sequence $b$ defined by $b_n = \norm{\varphi_n}$ and let the
  operator $B$ defined by $(B u)_n =  \langle \varphi_n, u\rangle$.

  Then, if $b$ is an $L^p$--admissible input element for $A$, $B$ is
  bounded from $U$ to $X_{-\theta}$ for all $\theta > \sfrac1{p'}$ and
  $B$ is an $L^p$--admissible control operator for $A$.
\end{proposition}
\begin{proof}
  From the proof of Theorem~\ref{thm:alpha-Carleson-Poisson}, we know
  that $b$ is an element of $X_{-\theta}$ for all
  $\theta>\sfrac1{p'}$. The elementary estimate $|\langle \varphi_n, u
  \rangle| \le \norm{u} \, \norm{\varphi_n}$ then implies that $B$ is
  linear and bounded from $U$ to $X_{-\theta}$. Now let $u \in L^p(0,
  \infty; U) = L^p(0, \infty, \ell_q)$ and let $u_j(\cdot)$ denote its
  coordinate functions. Then
\begin{align*}
    \biggnorm{ \int_0^t S(t{-}s) B u(s)\,ds }_{\ell_q} 
= & \;\biggl( \sum_{n=1}^\infty \biggl| \int_0^t e^{-\la_n (t{-}s)} \langle
      \varphi_n, u(s) \rangle \,ds \biggr|^q \biggr)^\sfrac1q \\
= & \;\biggl( \sum_{n=1}^\infty \biggl| \bigl\langle \varphi_n, \int_0^t e^{-\la_n (t{-}s)} 
       u(s) \bigr\rangle \,ds \biggr|^q \biggr)^\sfrac1q \\
\le& \;\biggl( \sum_{n=1}^\infty \norm{\varphi_n}^q \biggnorm{ \int_0^t e^{-\la_n (t{-}s)} u(s) \,ds}_U^q \biggr)^\sfrac1q \\
 = & \;\biggl( \sum_{n=1}^\infty \norm{\varphi_n}^q \sum_{j=1}^\infty
    \biggl| \int_0^t e^{-\la_n (t{-}s)} u_j(s) \,ds \biggr|^q \biggr)^\sfrac1q \\
 = & \;\biggl( \sum_{j=1}^\infty \biggl[ \sum_{n=1}^\infty  \biggl| \int_0^t
     e^{-\la_n (t{-}s)} \norm{\varphi_n} u_j(s) \,ds \biggr|^q \biggr]
     \biggr)^\sfrac1q.
\intertext{Now, by assumption, we may estimate}
\le& \; C \,\biggl( \sum_{j=1}^\infty \bignorm{ u_j }_{L^p}^q \biggr)^\sfrac1q
     = C \,\norm{ u }_{\ell_q(L^p(0,\infty))} \\
\le& \; C \,\norm{ u }_{L^p(0,\infty; U)} 
\end{align*}
where the last estimate (namely that $\norm{f}_{\ell_q(L^p)} \le
\norm{f}_{L^p(\ell_q)}$ if $p \le q$) is a easy corollary of Minkowski's
inequality, see e.g. \cite[Exercise VI.11.14]{DunfordSchwartz:Part1}.
\end{proof}

Notice that combining the Proposition with
Theorems~\ref{thm:alpha-Carleson-Poisson} and
\ref{thm:alpha-Carleson-direkt} yields a sufficient condition for
$L^p$--admissibility for control operators $B: U \to X_{-1}$ via a
Carleson measure criterion under suitable conditions on $p$ and $q$.

\section{$L^p$--admissibility of reciprocal systems}\label{sec:reciprocal}

Suppose that $\al=\sfrac{p'}q \le 1$ and $p\le2$.
Then, if $b$ is $L^{p'}$--admissible for $A$, 
Theorem~\ref{thm:alpha-Carleson-direkt} tell us that the discrete 
measure $\mu = \sum_{n} \bigl|\tfrac{b_n}{\la_n}\bigr|^q \delta_{\la_n^{-1}}$
is embedding $\al$--Carleson. 
Since the numbers $\tfrac{b_n}{\la_n}$ are the coefficients of $A^{-1}b$ and since 
the support of the measure is equal to the set of eigenvalues of $A^{-1}$, 
Theorem~\ref{thm:alpha-Carleson-Poisson} implies that $A^{-1}b$ is 
an $L^p$ admissible input element for $A^{-1}$. 

This observation can be formalised in the following question: given a
sectorial operator $A$ and a (control) operator $B \in \BOUNDED(U,
X_{-1})$, does $L^{p'}$--admissibility of the pair $A, B$ imply the
$L^p$--admissibility of the pair $A^{-1}, A^{-1}B$ ?  This question is
not new. In fact, the system
\begin{equation} \label{eq:reciprocal}
   z'(t) + A^{-1} z(t) = A^{-1}B u(t), \qquad z(0) = z_0.
\end{equation}
is called the {\em associated reciprocal linear system} to
(\ref{eq:control-system}).  By considering Lyapunov equations 
{\sc  Curtain} shows in the case of Hilbert spaces $X$ and $U$ that, under
the assumption that $0\in \varrho(A)$, an operator $B$ is
infinite-time $L^2$--admissible for $A$ if and only if $A^{-1}B$ is
infinite-time $L^2$--admissible for $A^{-1}$ (see \cite[Theorem
5]{Curtain:reciprocal}). We shall discuss generalisations of this
result below.

\bigskip

{\bf Analytic semigroups on Banach spaces.}  The famous Weiss
conjecture (see \cite{Weiss:conjectures}) states that
$L^2$--admissibility of a control operator $B$ between Hilbert spaces
$U$ and $X$ is characterised by the boundedness of the set
\begin{equation} \label{eq:weiss-cond}
  \bigl\{\, \Re(\la)^\einhalb (\la+A)^{-1} B \suchthat \Re(\la) >0 \,\bigr\}  
\end{equation}
in $\BOUNDED(U, X)$. The conjecture holds for normal semigroups and $U=\CC$
but fails even in case of $U=\CC$ for general semigroups on Hilbert
spaces (see \cite{JacobZwart:counterexable-Weiss-conj}).  For a survey
on positive and negative results concerning the Weiss conjecture we
refer to \cite{JacobPartington:survey}.

For analytic semigroups, $L^2$--admissibility is characterised by
(\ref{eq:weiss-cond}) even in Banach space under an inevitable
additional condition (see \cite{LeMerdy:weiss-conj}). This
characterisation has subsequently been generalised to $L^p$--norms
with certain weights (see \cite{HaakLeMerdy} and
\cite{HaakKunstmann:weighted-Lp-admiss}). We briefly summarise some
required notions and results: \\

\smallskip

For sectorial operators $A$ there exists a natural functional calculus
on the algebra of functions
\[
H_0^\infty(\Sec{\om}) = \bigl\{ f\in H^\infty(\Sec{\om}) \suchthat 
\exists c,s>0:\; \bigl|f(z)\bigr| \le c \,\min(|z|^s, |z|^{-s}) \bigr\}
\]
that is given by $f(A) := \int_\Ga f(z) R(z, A)\,dz$. Here, $\Ga$ is
the positively orientated boundary of a sector $S(\theta)$ containing
the spectrum of $A$ (see \cite{Haase:Buch,McIntosh:H-infty-calc} for
more details).  Since $B = A^\einhalb$ always satisfies condition
(\ref{eq:weiss-cond}), we find that whenever the Weiss conjecture
holds true, one has
\begin{equation}   \label{eq:square-function-estimate}
    \int_0^\infty \bignorm{ \varphi(tA) x }^2 \tfrac{dt}t 
=   \int_0^\infty \bignorm{ (tA)^\einhalb S(t) x }^2 \tfrac{dt}t 
\le m^2 \norm{x}^2
\end{equation}
for $\varphi(z) = z^\einhalb \exp(-z)$. By a result of {\sc
  M${}^c$Intosh} and {\sc Yagi} \cite[Theorem
5]{McIntoshYagi:without}, inequality
(\ref{eq:square-function-estimate}) does not depend on the particular
choice of the function $\varphi \in H_0^\infty(S(\om))\weg\{0\}$. This
result extends to Banach spaces and $L^p$--norms for all $p\in
[1,\infty]$. Such estimates are called $L^p_\ast$--estimates, where
the subscript refers to the measure $dt/t$ on $\RR_+$. They
automatically hold for $p=\infty$. Moreover, $L^p_\ast$ estimates for
$A$ imply $L^q_\ast$-estimates for all $q\ge p$ (see \cite[Remark
1.6]{HaakKunstmann:weighted-Lp-admiss}).  This is due to the fact that
$L^p_\ast$-estimates for $A$ are equivalent to $X \emb (\dot X_{-1},
\dot X_1))_{\einhalb, p}$ whereas $L^q_\ast$-estimates for $A'$ are
equivalent to $(\dot X_{-1}, \dot X_1)_{\einhalb, q'} \emb X$.  We
quote a characterisation of admissibility for control operators from
\cite{HaakKunstmann:weighted-Lp-admiss}. Recall that a sectorial
operator with dense range is actually injective (see \cite[Thm.
3.8]{McIntoshYagi:without}).

\begin{theorem}\label{thm:control-Lp-zul}
Let $p\in[1,\infty)$ and $A$ be a densely defined sectorial operator
of type $\om<\pihalbe$ with dense range on $X$. For $B \in
\BOUNDED(U,X_{-1})$ consider the set
\begin{equation}\label{eq:W_B}
  W_{A,B}^{(p)} := \bigl\{ \la^{\sfrac{1}{p}}\,(\la+A)^{-1}B: \; \la>0 \bigr\}  \subseteq \BOUNDED(U,X).
\end{equation}
Then the following assertions hold:
\begin{enumerate}
\item \label{item:control-Lp-zul-a} 
      If $B$ is $L^p$--admissible for $A$, then $W_{A,B}^{(p)}$ is bounded.
\item \label{item:control-Lp-zul-b} 
      If $W_{A,B}^{(p)}$ is bounded and the adjoint operator $A'$ satisfies
      $L^{p'}_{\ast}$--estimates then $B$ is $L^p$--admissible for $A$.
\end{enumerate}
\end{theorem}

As in the case $p=2$, the additional condition in
\ref{item:control-Lp-zul-b} is optimal.

\begin{corollary}\label{cor:reciprocal-systems}
Let $p\in[1,\infty)$ and $A$ be a densely defined sectorial operator
of type $\om<\pihalbe$ with dense range on $X$ and let $B \in
\BOUNDED(U,X_{-1})$. 
\begin{enumerate}
\item If $B$ is $L^p$--admissible for $A$ and if $A'$ satisfies 
      $L^{p}_{\ast}$--estimates on $X'$, then
      $A^{-1}B$ is $L^{p'}$--admissible for $A^{-1}$.
\item If $A^{-1}B$ is $L^{p'}$--admissible for $A^{-1}$ and if 
      $A'$ satisfies $L^{p'}_{\ast}$--estimates on $X'$, then
      $B$ is $L^p$--admissible for $A$.
\end{enumerate}
\end{corollary}
\begin{proof}
  Notice that whenever $\varphi \in H_0^\infty(S(\theta))$, the
  function $\psi(z) = \varphi(\sfrac1z) \in H_0^\infty(S(\theta))$.
  By functional calculus, $\varphi(tA) = \psi(t^{-1} A^{-1})$ (see
  \cite[Proposition 2.4.1]{Haase:Buch}) and thus $A$ admits
  $L^r_\ast$--estimates if, and only if its inverse $A^{-1}$ does.  In
  view of Theorem~\ref{thm:control-Lp-zul} it is therefore sufficient
  to consider the expressions
\[
\la^{\sfrac1p} (\la+A)^{-1}B 
= \la^{\sfrac1p-1} (\la^{-1}+A^{-1})^{-1} A^{-1}B 
= z^{\sfrac1{p'}} (z+A^{-1})^{-1} A^{-1}B 
\]
where $\la^{-1}=z$. Now, if $B$ is $L^p$ admissible, $W_{A,B}^{(p)}$
is bounded, whence $W_{A^{-1},A^{-1} B}^{(p')}$ is bounded. If $A'$
satisfies $L^p_\ast$--estimates, then $A^{-1} B$ is
$L^{p'}$--admissible for $A$.  For the converse we need to impose
$L^{p'}_\ast$--estimates on $A'$ by the same argument.
\end{proof}

\begin{remark}
  We can use the corollary to compare
  Theorems~\ref{thm:alpha-Carleson-Poisson} and
  \ref{thm:alpha-Carleson-direkt} in the following case: let $A$
  generate a bounded analytic semigroup and assume that
\[
\mu = \sum_{n=1}^\infty \delta_{\la_n} |b_n|^q 
= \sum_{n=1}^\infty \delta_{(\la_n^{-1})^{-1}} \Bigl| \frac{\la_n^{-1} b_n}{\la_n^{-1}}\Bigr|^q 
\]
is embedding $\al$--Carleson with $\al=\sfrac{p'}q$, so that
Theorem~\ref{thm:alpha-Carleson-direkt} yields ${L^p}'$--admissibility
of $A^{-1}b$ for $A^{-1}$. If $p\le q$
Corollary~\ref{cor:reciprocal-systems} implies that $b$ is
$L^p$--admissible for $A$ since in this case the Lorentz space
$\ell_{q,p}$ embeds into $\ell_{q,q} = \ell_q$. Both theorems coincide
therefore in this case and it is remarkable, that the restriction
$p\le 2$ in Theorem~\ref{thm:alpha-Carleson-Poisson} may be weakened
to $p\le q$ if $q>2$, thus allowing in particular the natural choice
$p=q$ for all $p\in(1,\infty)$. In view of
Example~\ref{ex:examples}~\ref{item:ex_c} analyticity of the semigroup
seems to be a necessary restriction to guarantee this improvement of
Theorem~\ref{thm:alpha-Carleson-Poisson}.

If, conversely, $b$ is $L^p$--admissible for $A$,
Corollary~\ref{cor:reciprocal-systems} guarantees the
${L^p}'$--admissibility of $A^{-1}b$ for $A^{-1}$ when $p' \le q$ and
Theorem~\ref{thm:alpha-Carleson-direkt} yields that $\mu$ is geometric
$\al$--Carleson, so both theorems again coincide.
\end{remark}

\bigskip

{\bf Arbitrary semigroups on Banach spaces.}  In the sequel we want to
get rid of the assumption of analyticity and give a more direct
argument for passing from $L^p$--admissibility of $B$ to
$L^{p'}$--admissibility of $A^{-1}B$ for $A^{-1}$. To this end we
assume that $S(t) = \exp(-tA)$ is uniformly exponentially stable.
Then, by the {\em Phillips functional calculus} (see e.g.
\cite{Haase:Buch}) one obtains different representation formulas for
$\exp(-tA^{-1})$ that may turn out to be useful for different
problems. The most interesting (to us) are
\begin{equation}\label{eq:darstellung-zwart}
  \exp(-tA^{-1})x 
= x - \int_0^\infty \left(\tfrac{t}s\right)^\einhalb J_1(2\sqrt{st}) \; S(s) x\,ds
\end{equation}
for $x\in X$ where $J_1$ is the first Bessel function of first order.
This formula was given in
\cite{GomilkoZwart,GomilkoTomilovZwart,Zwart:inverse,Zwart:growth}.
It is an immediate consequence of \cite[Formula
(5.67)]{Oberhettinger:Laplace}.  Moreover we have
\begin{equation}\label{eq:darstellung-einhalb}
  A^{-\einhalb} \exp(-tA^{-1})x 
=  \int_0^\infty \tfrac1{\sqrt{\pi s}} \cos(2\sqrt{st})\;S(s) x\,ds
\end{equation}
\begin{equation}\label{eq:darstellung-A-potenz}
  A^{-\nu-1} \exp(-tA^{-1})x 
= \int_0^\infty \left(\tfrac{t}s\right)^{\sfrac{\nu}2} I_\nu(2\sqrt{st}) 
  \; S(s) x\,ds, \qquad \nu>-1
\end{equation}
for $x\in X$ as a consequence of 
\cite[Formulas (5.70) and (5.75)]{Oberhettinger:Laplace}. Here, $I_\nu$ 
is the modified Bessel function of the first kind, of order $\nu$ (see
e.g. \cite{AbramowitzStegun,Watson:Bessel} for more details on Bessel 
functions).

Since for general exponentially stable semigroups, the decay 
rate 
\[
   \norm{  \exp(-tA^{-1})A^{-1}x } = O((1+t)^{-\sfrac14}) \norm{x}
\]
is optimal (see \cite[Theorem 3.3 and Example 3.5]{Zwart:inverse}, one
cannot expect more than that $L^p$--admissibility of $B$ implies
finite-time $L^{p'}$--admissibility of $A^{-1}B$ for $A^{-1}$. In this
sense the following result is optimal.

\begin{theorem}
  Let $-A$ generate an exponentially stable semigroup $S(t)_{t\ge0}$
  on a Banach space $X$, suppose that $A$ has dense range and that $B
  \in \BOUNDED(U, X_{-1})$ is $L^p$--admissible for $A$, $p\in[1,\infty]$.
  Then $A^{-1}B$ is finite-time $L^{p'}$--admissible for $A^{-1}$.
\end{theorem}
\begin{proof}
  Let $T>0$. Then we may find $\eps>0$ such that $B$ is
  $L^p$--admissible for the shifted semigroup $e^{+\eps \,t}S(t)$ as
  well. By formula (\ref{eq:darstellung-A-potenz}) with $\nu=0$ we
  obtain
\begin{eqnarray*}
&   & \Bignorm{ \int_0^T \exp( -sA^{-1} ) A^{-1}B u(s)\,ds}_X \\
& = & \Bignorm{ \int_0^T \Bigl(\int_0^\infty I_0(2\sqrt{st}) S(t)\,dt\Bigr) 
       B u(s)\,ds }_X \\
& = & \Bignorm{ \int_0^\infty e^{+\eps \,t} S(t) B \Bigl(
        e^{-\eps \,t} \int_0^T I_0(2\sqrt{st}) u(s)\,ds \Bigr) \,dt }_X \\
&\le& M \Bignorm{e^{-\eps \,t} \int_0^T I_0(2\sqrt{st}) 
        u(s)\,ds}_{L^p((0,\infty), U)} \\
&\le& M \Bignorm{e^{-\eps \,t} \int_0^T \bigl|I_0(2\sqrt{st})\bigr| 
         \,\norm{ u(s)}\,ds }_{L^p(0,\infty)} \\
&\le& M \norm{u}_{L^{p'}(0,T)}  \Bignorm{ t\mapsto  e^{-\eps \,t} \bignorm{
          s\mapsto I_0(2\sqrt{st}) }_{L^{p}(0,T)} }_{L^p(0,\infty)}
   =: M C_T \; \norm{u}_{L^{p'}(0,T)}.
\end{eqnarray*}
From \cite[Formula (9.6.16)]{AbramowitzStegun}, we deduce the
following integral representation
\[
I_0(x) = \tfrac1\pi \int_0^\pi e^{x\, \cos(\theta)} \,d\theta
       = \tfrac2\pi e^x \int_0^1 \exp(-2u^2 x) (1-u^2)^{-\einhalb}\,du
\]
by letting $u = \sin(\sfrac\theta2)$ (see also \cite[p.
204]{Watson:Bessel}).  In particular, $|I_0(x)|\le e^x$, $x>0$. In the
case $p=\infty$ boundedness of $C_T$ is immediate. In case the
$p<\infty$,
\[
f(t) := \bignorm{s\mapsto I_0(2\sqrt{st}) }_{L^{p}(0,T)}^p 
   \le \tfrac2{tp^2} (1+e^{p\sqrt{tT}}(p\sqrt{tT}-1)).
\]
Now, since $\lim_{t\to 0} f(t) = T$ it is clear that $e^{-\eps
  \,t}f(t)$ is integrable on $(0,\infty)$ for any $\eps>0$. Thus $C_T
< \infty$ and so $A^{-1}B$ is finite-time $L^{p'}$--admissible for
$A^{-1}$.
\end{proof}

\subsection*{Acknowledgement}
 We are indepted to Birgit Jacob, Peer Christian Kunstmann and
 Elizabeth Strouse for helpful remarks and discussions.

\def\SUBMITTED{Submitted}
\def\TOAPPEAR{To appear in }
\def\PREPARATION{In preparation }

\def\cprime{$'$}
\providecommand{\bysame}{\leavevmode\hbox to3em{\hrulefill}\thinspace}

\end{document}